\documentclass[12pt,a4paper]{article}
\usepackage{epsfig}
\usepackage{amsmath}
\usepackage{latexsym}
\usepackage{amssymb}
\usepackage{paralist}
\usepackage{amsfonts}
\newtheorem{thm}{Theorem}[section]
\newtheorem{definition}[thm]{Definition}
\newtheorem{prop}[thm]{Proposition}
\newtheorem{lem}[thm]{Lemma}
\newtheorem{cor}[thm]{Corollary}

\newcommand{\pf}{\noindent {\bf Proof: \ }}
\newcommand\Q[1]{$#1$ \quad}

\font\msbm=msbm10 at 12pt

\newcommand{\FF}{\mbox{\msbm F}}

\begin{document} 
\title{The lengths \\ of\\ Hermitian Self-Dual\\ Extended Duadic Codes}

\author{Lilibeth Dicuangco\footnote{Mathematics Department, University of the Philippines, Diliman, Quezon 
City, 1101 Philippines, {\tt ldicuangco@math.upd.edu.ph}}, \: \setcounter{footnote}{6} Pieter 
Moree\footnote{Max-Planck-Institut, Vivatsgasse 7, D-53111 Bonn, Germany, {\tt moree@mpim-bonn.mpg.de}},  
\: \setcounter{footnote}{3} Patrick Sol\'e\footnote{CNRS, I3S ESSI, BP 145, Route des Colles, 06 903 
Sophia Antipolis, France, {\tt sole@essi.fr}}}

\date{} 
\maketitle 
\date{} 
\begin{abstract}

\noindent Duadic codes are a class of cyclic codes that generalizes quadratic residue codes from prime to 
composite lengths. For every prime power $q$, we characterize integers $n$ such that there is a duadic 
code of length $n$ over $\FF_{q^{2}}$ with an Hermitian self-dual parity-check extension. We derive 
asymptotic estimates for the number of such $n$ as well as for the number of lengths for which duadic 
codes exist.
\end{abstract}

\noindent {\bf Keywords}: Duadic codes, Extended cyclic codes, Hermitian self-dual codes, Splittings.

\noindent {\bf Mathematics Subject Classification}: 11N64, 94B15, 11N37

\section{Introduction}

Duadic codes are a family of cyclic codes over fields that generalize quadratic residue codes to composite 
lengths. For a general introduction, see \cite{fundamentals}, \cite{leon et al} and \cite{smid}. It can be 
determined when an extended duadic code is self-dual for the Euclidean scalar product 
(\cite{fundamentals}). In this work, we study for which $n$ there exist duadic codes over $\FF_{q^{2}}$ of 
length $n$ the extension of which by a suitable parity-check is self-dual for the Hermitian scalar product 
$\sum_{i=1}^{n+1} x_{i}y_{i}^{q}$.

First, we characterize the Hermitian self-orthogonal cyclic codes by their defining sets (Proposition 
\ref{Thm 4.4.11.1}), then the duadic codes (Proposition \ref{Thm 6.4.1.1}). Next, we study under what 
conditions the extension by a parity-check of a duadic code is Hermitian self-dual (Proposition \ref{Thm 
6.4.14.1}). Finally, we derive by elementary means an arithmetic condition bearing on the divisors of $n$ 
(Theorem \ref{Thm 3.17}) for the previous situation. This condition was arrived at in \cite{conchita} 
using representation theory of groups. In an appendix, we derive asymptotic estimates for $x$ large on 
$A_{q}(x)$, the number of integers $\leq x$ that are split by the multiplier $\mu_{-q}$, and on 
$D_{q}(x)$, the number of possible lengths $\leq x$ of a duadic code. The proofs are based on analytic 
number theory.

\section{Preliminaries}

We assume the reader is familiar with the theory of cyclic codes (see e.g., \cite{intro in handbook}, 
\cite{fundamentals}). Let $q$ be a power of a prime $p$ and let $\FF_{q}$ denote the Galois field with $q$ 
elements. Let $n$ be a positive integer such that gcd$(n,q)=1$.  Let $\mathcal{R}_n = \FF_{q}\lbrack x 
\rbrack / (x^{n}-1)$. We view a cyclic code over $\FF_{q}$ of length $n$ as an ideal in $\mathcal{R}_n$. 

Let $0 < s < n$ be a nonnegative integer. Let $ C_{s} = \lbrace s, sq, sq^{2}, \ldots, sq^{r_{s}-1} 
\rbrace$ , where $r_{s}$ is the smallest positive integer such that $sq^{r_{s}} \equiv s \;(\mbox{mod } 
n)$. The coset $C_{s}$ is called the \emph{$q$-cyclotomic coset of $s$ modulo $n$}. The subscript of 
$C_{s}$ is usually taken to be the smallest number in the set and is also taken as the coset 
representative. The distinct $q$-cyclotomic cosets modulo $n$ partition the set $\lbrace 0,1,2,\ldots,n-1 
\rbrace$.

Let $\alpha$ be a primitive $n$th root of unity in some extension field of $\FF_{q}$. A set $T \subseteq 
\lbrace 0, 1, 2, \ldots, n-1 \rbrace$ is called the \emph{defining set} (relative to $\alpha$) of a cyclic 
code $C$ whenever $c(x) \in C$ if and only if $ c({\alpha^{i}}) = 0$ for all $i \in T$. In this paper, we 
assume implicitly that an $n$th root of unity has been fixed when talking of defining sets.

A ring element $e$ such that $e^{2}=e$ is called an idempotent. Since gcd$(n,q)=1$, the ring 
$\mathcal{R}_n$ is semi-simple. Thus, by invoking the Wedderburn Structure theorems, we can say that each 
cyclic code in $\mathcal{R}_n$ contains a unique idempotent element which generates the ideal. 
Alternatively, this fact has also been proven directly in \cite[Theorem 4.3.2]{fundamentals}. We call this 
idempotent element the \emph{generating idempotent} (or \emph{idempotent generator}) of the cyclic code.

Let $a$ be an integer such that gcd$(a,n)=1$. We define the function $\mu_{a}$, called a 
\emph{multiplier}, on $\lbrace 0,1,2, \cdots, n-1 \rbrace$ by $i\mu_{a}\equiv ia \;(\mbox{mod } n)$. 
Clearly,  $\mu_a$ gives a permutation of the coordinate positions of a cyclic code of length $n$. Note 
that this is equivalent to the action of $\mu_{a}$ on $\mathcal{R}_n$ by $f(x)\mu_{a} \equiv f(x^{a}) 
\;(\mbox{mod } x^{n}-1)$.

If $C$ is a code of length $n$ over $\FF_{q}$, we define a \emph{complement of C} as a code $C^{c}$ such 
that $C+C^{c}=\FF_{q}^{n}$ and $C \cap C^{c} = \lbrace \mathbf{0} \rbrace$. In general, a complement of a 
code is not unique. But it is easy to show that if $C$ is cyclic, then $C^{c}$ is unique and that it is 
also cyclic (see e.g., Exercise 243, \cite{fundamentals}). In this case, we call $C^{c}$ \emph{the cyclic 
complement of $C$}.

\section{Cyclic Codes over $\FF_{q^{2}}$}

We now consider cyclic  codes over the Galois field $\FF_{q^{2}}$, where $q$ is a power of a prime $p$. In 
this case, we note that $\mathcal{R}_n = \FF_{q^{2}}\lbrack x \rbrack / (x^{n}-1)$.

\subsection{Idempotents in $\mathcal{R}_{n}$}

Consider the involution $\; \bar{   }: z \mapsto z^{q}$ defined on $\FF_{q^{2}}$. We extend this map 
component-wise to $\FF_{q^{2}}^{n}$. For an element $a(x)=a_{0} + a_{1}x + \cdots + a_{n-1}x^{n-1}$ in 
$\mathcal{R}_{n}$, we set $\overline{a(x)} = a_{0}^{q} + a_{1}^{q}x + \cdots + a_{n-1}^{q}x^{n-1}$.

Let $C$ be a code of length $n$ over $\FF_{q^{2}}$. We define the \emph{conjugate of $C$} to be the code 
$\overline{C} = \lbrace \overline{\mathbf{c}} \mid \mathbf{c} \in C \rbrace$. It can easily be shown that 
if $C$ is a cyclic code with generating idempotent $e(x)$, then $\overline{C}$ is also cyclic and its 
generating idempotent is $\overline{e(x)}$.

Suppose we list all the distinct $q^{2}$-cyclotomic cosets modulo $n$ in the following way:
\begin{equation*}
C_{1}, C_{2},\ldots,C_{k},D_{1},D_{2},\ldots,D_{l},E_{1},E_{2},\ldots,E_{l},
\end{equation*}
such that
\begin{eqnarray*}
 			C_{i} = qC_{i}	\; \;   \mbox{ for } 1 \leq i \leq k & \; \; \mbox{and} \; \; & E_{i} = qD_{i} 
\; \;  \mbox{ for } 1 \leq i \leq l .
\end{eqnarray*}
By Corollary  4.3.15 of \cite{fundamentals}, an idempotent in $\mathcal{R}_{n}$ has the form
\begin{equation}
e(x) = \sum_{j=1}^{k} a_{j} \sum_{i \in C_{j}}x^{i} + \sum_{j=1}^{l} b_{j} \sum_{i \in D_{j}}x^{i} + 
\sum_{j=1}^{l} c_{j} \sum_{i \in E_{j}}x^{i}.
\end{equation}
Thus,
\begin{eqnarray*}
e(x) & = & e(x)^{q} \\
     & = & \sum_{j=1}^{k} a_{j}^{q} \sum_{i \in C_{j}}x^{qi} + \sum_{j=1}^{l} b_{j}^{q} \sum_{i \in 
D_{j}}x^{qi} + \sum_{j=1}^{l} c_{j}^{q} \sum_{i \in E_{j}}x^{qi} \\
	& = & \sum_{j=1}^{k} a_{j}^{q} \sum_{i \in C_{j}}x^{i} + \sum_{j=1}^{l} b_{j}^{q} \sum_{i \in 
E_{j}}x^{i} + \sum_{j=1}^{l} c_{j}^{q} \sum_{i \in D_{j}}x^{i}.
\end{eqnarray*}
Hence,
\begin{eqnarray*}
a_{j}^{q} = a_{j} & 1 \leq j \leq k \: ; \\
b_{j}^{q} = c_{j} & 1 \leq j \leq l \: ,
\end{eqnarray*}
which implies
\begin{equation}
e(x) = \sum_{j=1}^{k} a_{j} \sum_{i \in C_{j}} x^{i} + \sum_{j=1}^{l}b_{j}\sum_{i \in D_{j}} x^{i} + 
\sum_{j=1}^{l}b_{j}^{q} \sum_{i \in D_{j}}x^{qi}.
\label{star1}
\end{equation}
Thus,
\begin{eqnarray*}
\overline{e(x)} & = & \sum_{j=1}^{k} a_{j}^{q} \sum_{i \in C_{j}} x^{i} + \sum_{j=1}^{l}b_{j}^{q}\sum_{i 
\in D_{j}} x^{i} + \sum_{j=1}^{l}b_{j}^{q^{2}} \sum_{i \in D_{j}}x^{qi} \\
           & = & \sum_{j=1}^{k} a_{j} \sum_{i \in C_{j}} x^{qi} + \sum_{j=1}^{l}c_{j}\sum_{i \in E_{j}} 
x^{qi} + \sum_{j=1}^{l}b_{j} \sum_{i \in D_{j}}x^{qi} \\
		   & = & e(x)\mu_{q}.
\end{eqnarray*}
This gives $\overline{C} = \langle \overline{e(x)} \rangle = \langle e(x)\mu_{q} \rangle = C\mu_{q}$ by 
Theorem 4.3.13 of \cite{fundamentals}.

\vskip .1in
The discussion above is summarized in the following proposition.

\begin{prop}
Let $C$ be a cyclic code over $\FF_{q^{2}}$ with generating idempotent $e(x)$. The following hold:
\begin{enumerate}
\item $e(x)$ has the form given in {\rm (\ref{star1})}.
\item $\overline{C}$ is cyclic with generating idempotent $\overline{e(x)}$.
\item $\overline{e(x)} = e(x)\mu_{q}$.
\item $\overline{C} = C\mu_{q}$.
\end{enumerate} 
\label{Thm 4.3.16.1}
\end{prop}

\subsection{Euclidean and Hermitian Duals}

Let $\mathbf{x}=(x_{0}, x_{1},\ldots,x_{n-1})$ and $\mathbf{y}=(y_{0},y_{1},\ldots,y_{n-1})$ be any 
vectors in $\FF_{q^{2}}^{n}$. Consider the involution $\bar{ }: z \mapsto z^{q}$ defined on $\FF_{q^{2}}$.
The \emph{Hermitian scalar product} of $\mathbf{x}$ and $\mathbf{y}$ is given by $\mathbf{x} \cdot 
\overline{\mathbf{y}} = \sum_{i=0}^{n-1} x_{i}\overline{y_{i}}$. If $C$ is a linear code over 
$\FF_{q^{2}}$, the \emph{Euclidean dual of $C$} is denoted $C^{\bot_{E}}$. The \emph{Hermitian dual of 
$C$} is $C^{\bot_{H}} = \lbrace \mathbf{u} \in \FF_{q^{2}}^{n} \mid \mathbf{u} \cdot \overline{\mathbf{w}} 
= 0 \mbox{ for all } \mathbf{w} \in C \rbrace$. We say that a code $C$ is \emph{Euclidean self-orthogonal} 
if $C \subseteq C^{\bot_{E}}$, and that  $C$ is \emph{Euclidean self-dual} if $C = C^{\bot_{E}}$. 
Similarly, $C$ is said to be \emph{Hermitian self-orthogonal} if $C \subseteq C^{\bot_{H}}$, and that $C$ 
is \emph{Hermitian self-dual} if $C = C^{\bot_{H}}$.

Let $f(x)=f_{0} + f_{1}x + \cdots + f_{r}x^{r} \in \FF_{q^{2}} \lbrack x \rbrack $. The \emph{reciprocal 
polynomial} of $f(x)$ is the polynomial $f^{*}(x) = x^{r}f(x^{-1}) = x^{r}(f(x)\mu_{-1}) = f_{r} + 
f_{r-1}x + \cdots + f_{0}x^{r}$.

\begin{lem}
Let $\mathbf{a}= (a_{0}, a_{1},\ldots, a_{n-1})$, $\mathbf{b}= (b_{0}, b_{1},\ldots, b_{n-1})$ be vectors 
in $\FF_{q^{2}}^{n}$ with associated polynomials $a(x)$ and $b(x)$. Then $\mathbf{a}$ is Hermitian 
orthogonal (similarly Euclidean orthogonal) to $\mathbf{b}$ and all its cyclic shifts if and only if 
$a(x)\overline{b^{*}(x)} = 0$ (similarly $a(x)b^{*}(x)=0$) in $\mathcal{R}_{n}$.
\label{Lem 4.4.8.1}
\end{lem}

\pf  See Lemma 4.4.8 of \cite{fundamentals} for the Euclidean case. The proof for the case of Hermitian 
orthogonality follows a similar argument and is omitted. \hfill \Q{\Box}
\vskip .1in
\noindent Recall that a vector $\mathbf{a} = (a_{0}, a_{1}, \ldots, a_{n-1})$ of $\FF_{q^{2}}^{n}$ is 
called an \emph{even-like vector} if $\sum_{i=0}^{n-1} a_{i} = 0$. A code $C$ is called an \emph{even-like 
code} if all its codewords are even-like, otherwise it is called \emph{odd-like}.
The following lemma appears as Exercise 238 of \cite{fundamentals} and its proof is left to the reader.

\begin{lem}
Let $C$ be a cyclic code over $\FF_{q^{2}}$ with defining set $T$ and generator polynomial $g(x)$. Let 
$C_{e}$ be the subcode of $C$ consisting of all the even-like vectors in $C$. Then:
\begin{enumerate}
\item $C_{e}$ is cyclic and has defining set $T \cup \lbrace 0 \rbrace$.
\item $C = C_{e}$ if and only if $0 \in T$ if and only if  $g(1)=0$.
\item If $C \neq C_{e}$, then the generator polynomial of $C_{e}$ is $(x-1)g(x)$.
\end{enumerate}
\label{exercise 238}
\end{lem}

\noindent The following propositions generalize some results on Euclidean duals of cyclic codes over an 
arbitrary finite field to Hermitian duals of cyclic codes over $\FF_{q^{2}}$.

\begin{prop}
Let $C$ be a cyclic code of length $n$ over $\FF_{q^{2}}$ with generating idempotent $e(x)$ and defining 
set $T$. The following hold:
\begin{enumerate}
\item $C^{\bot_{H}}$ is a cyclic code and $C^{\bot_{H}} = C^{c}\mu_{-q}$.
\item $C^{\bot_{H}}$ has generating idempotent $1 - e(x)\mu_{-q}$.
\item If $\mathcal{N} = \lbrace 0,1,2, \ldots , n-1 \rbrace $, then $\mathcal{N} \setminus (-q)T$ mod $n$ 
is the defining set for $C^{\bot_{H}}$.
\item Precisely one of $C$ and $C^{\bot_{H}}$ is odd-like and the other is even-like.
\end{enumerate}
\label{Thm 4.4.9.1}
\end{prop}

\pf  Let $\mathbf{a}= (a_{0}, a_{1},\ldots, a_{n-1}) \in C$. Denote by $\mathbf{a}^{(i)}$ the $i^{th}$ 
cyclic shift of $\mathbf{a}$. 
By assumption $\mathbf{a}^{(i)} \in C$ for all $i$.
Let $\mathbf{b}= ( b_{0}, b_{1},\ldots, b_{n-1} ) \in C^{\bot_{H}}$.
For $i = 0,1,\ldots,n-1$, we have  $\mathbf{b}^{(i)} \cdot \overline{\mathbf{a}} = \mathbf{b} \cdot 
\overline{\mathbf{a}^{n-i}} = 0$. 
Thus $C^{\bot_{H}}$ 
 is cyclic.
Note that $C^{\bot_{H}} = \overline{C}^{\bot_{E}}$ and $\overline{C}^{\bot_{E}}=\overline{C}^{c}\mu_{-1}$ 
(Theorem 4.4.9 of \cite{fundamentals}). 
It can easily be shown that $\overline{C}^{c} = \overline{C^{c}}$ and so using Proposition \ref{Thm 
4.3.16.1} we have $\overline{C}^{c} = C^{c}\mu_{q}$.
Hence $C^{\bot_{H}}=\overline{C}^{c}\mu_{-1} = C^{c}\mu_{q}\mu_{-1} = C^{c}\mu_{-q}$, proving  part 1.

\vskip 0in
Using Theorem 4.4.6 and Theorem 4.3.13 of \cite{fundamentals}, the generating idempotent for $C^{\bot_{H}} 
= C^{c}\mu_{-q}$ is $(1-e(x))\mu_{-q} = 1 - e(x) \mu_{-q}$. Thus part 2 holds.

\vskip 0in
The defining set for $C^{c}$ is $\mathcal{N} \setminus T$ (Theorem 4.4.6 of \cite{fundamentals}) and hence 
applying Corollary 4.4.5 of \cite{fundamentals} the defining set for $C^{\bot_{H}}$ is 
$(-q)^{-1}(\mathcal{N} \setminus T)  = \mathcal{N} \setminus (-q)^{-1}T \mbox{ mod } n$. 
Since $\mu_{-q}^{2} = \mu_{(-q)^{2}} = \mu_{q^{2}}$ 
fixes each $q^{2}$-cyclotomic coset and $(-q)^{-1}T$ is a union of $q^{2}$-cyclotomic coset (using Theorem 
4.4.2 of \cite{fundamentals}), it follows that $(-q)^{-1}T = (-q)^{2}(-q)^{-1}T = (-q)T \mbox{ mod } n$. 
Thus the defining set for $C^{\bot_{H}}$ is $\mathcal{N} \setminus (-q)T \mbox{ mod } n$. This proves part 
3.

\vskip 0in
Lastly, since exactly one of $T$ and $\mathcal{N} \setminus (-q)T$ contains $0$, part 4 follows from part 
3 and Lemma \ref{exercise 238} .  \hfill \Q{\Box}

\vskip .1in
\noindent The following lemma is from Exercise 239 of \cite{fundamentals}.

\begin{lem}
Let $C_{i}$ be a cyclic code of length $n$ over $\FF_{q^{2}}$ with defining sets $T_{i}$ for $i=1,2$. 
Then:
\begin{enumerate}
\item $C_{1} \cap C_{2}$ has defining set $T_{1} \cup T_{2}$.
\item $C_{1} + C_{2}$ has defining set $T_{1} \cap T_{2}$.
\item $C_{1} \subseteq C_{2} \iff T_{2} \subseteq T_{1}$.
\end{enumerate}
\label{Exercise 239}
\end{lem}

\begin{prop}
Let $C$ be a Hermitian self-orthogonal cyclic code over $\FF_{q^{2}}$ of length $n$ with defining set $T$.
Let $ \; C_{1}, C_{2},\ldots,C_{k},\; D_{1},D_{2},\ldots,D_{l}, \; E_{1},E_{2},\ldots,E_{l} \;$ be all the 
distinct  $q^{2}$-cyclotomic cosets modulo $n$ partitioned such that $C_{i} = C_{i}\mu_{-q}$ for  $1 \leq 
i \leq k$ and $D_{i} = E_{i}\mu_{-q}$ for  $1 \leq i \leq l$.
Then the following hold:
\begin{enumerate}
\item $C_{i} \subseteq T \mbox{ for } 1 \leq i \leq k$, and at least one of $D_{i}$ or $E_{i}$ is 
contained in $T$ for each $1 \leq i \leq l$.
\item $C$ is even-like.
\item $C \cap C\mu_{-q} = \lbrace 0 \rbrace$.
\end{enumerate}
Conversely, if $C$ is a cyclic code with defining set $T$ that satisfies part 1, then $C$ is an Hermitian 
self-orthogonal code.
\label{Thm 4.4.11.1}
\end{prop}

\pf Let $\mathcal{N}= \lbrace 0,1,2,\ldots,n-1 \rbrace $. 
By Proposition \ref{Thm 4.4.9.1}, $T^{\bot} = \mathcal{N} \setminus (-q)T \mbox{ mod } n$ is the defining 
set for $C^{\bot_{H}}$.
By assumption, $C \subseteq C^{\bot_{H}}$ and so $\mathcal{N} \setminus (-q)T \subseteq T$ by Lemma 
\ref{Exercise 239}.
If $C_{i} \not\subseteq T$ for some $i$, then $C_{i}\mu_{-q} \not\subseteq (-q)T$.
Since $C_{i} = C_{i}\mu_{-q}$, it follows that $C_{i} \subseteq \mathcal{N} \setminus (-q)T \subseteq T$, 
a contradiction. Thus $C_{i} \subseteq T$ for all $i$.
If $D_{i} \not\subseteq T$, then $E_{i} = D_{i}\mu_{-q} \not\subseteq (-q)T \mbox{ mod } n$.
Thus $E_{i} \subseteq \mathcal{N} \setminus (-q)T \subseteq T$. Hence part 1 holds.
\vskip 0in
To prove part 2, note that $\lbrace 0 \rbrace = C_{i}$ for some $i$. Hence $0 \in T$ by part 1. By Lemma 
\ref{exercise 238}, $C$ is even-like.
\vskip 0in
As noted in the proof of Proposition \ref{Thm 4.4.9.1}, $(-q)^{-1}T = (-q)T \mbox{ mod } n$. Using 
Corollary 4.4.5 of \cite{fundamentals}, $C\mu_{-q}$ has defining set $(-q)^{-1}T = (-q)T$.
Since $\mathcal{N} \setminus (-q)T \subseteq T$, it follows that $T \cup (-q)T = \mathcal{N}$. By Lemma 
\ref{Exercise 239}, $T \cup (-q)T$ is the defining set for $C \cap C\mu_{-q}$. Thus $C \cap C\mu_{-q} = 
\lbrace 0 \rbrace$, which proves part 3.
\vskip 0in
For the converse, assume $T$ satisfies part 1. We will show that $T^{\bot} \subseteq T$ which will imply 
that $C$ is Hermitian self-orthogonal.
By Proposition \ref{Thm 4.4.9.1}, $T^{\bot} = \mathcal{N} \setminus (-q)T \mbox{ mod } n$.
Note that $C_{i} \subseteq T \Longrightarrow C_{i} = C_{i}\mu_{-q} \subseteq (-q)T \Longrightarrow C_{i} 
\not\subseteq T^{\bot}$.
Hence $T^{\bot}$ is a union of some $E_{i}$'s and $D_{i}$'s.
If $D_{i} \subseteq T^{\bot} = \mathcal{N} \setminus (-q)T$, then $D_{i} \not\subseteq (-q)T \mbox{ mod } 
n$, implying that $(-q)D_{i} \not\subseteq T$. Since $(-q)D_{i}=E_{i}$, it follows that $E_{i} 
\not\subseteq T$. By part 1, $D_{i} \subseteq T$.
By a similar argument, it can be shown that if $E_{i} \subseteq T^{\bot}$, then $E_{i} \subseteq T$.
\hfill \Q{\Box}

\vskip .1in

\section{Duadic Codes}

Let $n$ be an odd positive integer. We let $\overline{j}(x) = \frac{1}{n}(1 + x + x^{2} + \cdots + 
x^{n-1})$, the generating idempotent for the repetition code of length $n$ over $\FF_{q}$. 

We first define duadic codes over arbitrary finite fields. Then we proceed to examine duadic codes over 
finite fields of square order. The goal of this section is to present some results concerning Hermitian 
orthogonality of duadic codes over such finite fields.

\subsection{Definitions and Basic Properties}

\begin{definition}
Let $e_{1}(x)$ and $e_{2}(x)$ be a pair of even-like idempotents and let $C_{1} = \langle e_{1}(x) \rangle 
$ and $C_{2} = \langle e_{2}(x) \rangle $. The codes $C_{1}$ and $C_{2}$ form a pair of \emph{{\bf 
even-like duadic codes}} if the following properties are satisfied: 

a) the idempotents satisfy $e_{1}(x) + e_{2}(x) = 1 - \overline{j}(x)$,

b) there is a multiplier $\mu_{a}$ such that $C_{1}\mu_{a} = C_{2} \mbox{ and } C_{2}\mu_{a} = C_{1}$.

\noindent To the pair of even-like codes $C_{1}$ and $C_{2}$, we associate a pair of \emph{{\bf odd-like 
duadic codes}}
$D_{1} = \langle 1 - e_{2}(x) \rangle$ and   $D_{2} = \langle 1 - e_{1}(x) \rangle$.
We say that the multiplier $\mu_{a}$ gives a \emph{{\bf splitting}} for the even-like duadic codes or for 
the odd-like duadic codes.
\label{duadicdef}
\end{definition}

\begin{thm}
{\rm (\cite{fundamentals}.)}
Let $C_{1}$ and $C_{2}$ be cyclic codes over $\FF_{q}$ with defining sets $T_{1} = \lbrace 0 \rbrace \cup 
S_{1}$ and $T_{2} = \lbrace 0 \rbrace \cup S_{2}$, respectively, where $0 \not\in S_{1}$ and $0 \not\in 
S_{2}$. Then $C_{1}$ and $C_{2}$ form a pair of even-like duadic codes if and only if the following 
conditions are satisfied: 

a) $S_{1}$ and $S_{2}$ satisfy $S_{1} \cup S_{2} = \lbrace 1,2, \ldots, n-1 \rbrace \mbox{ and } S_{1} 
\cap S_{2} = \emptyset$,

b) there is a multiplier $\mu_{b}$ such that $S_{1}\mu_{b} = S_{2} \mbox{ and } S_{2}\mu_{b} = S_{1}$ .
\label{Thm 6.1.5}
\end{thm}

\noindent If the conditions in the preceding theorem are satisfied, we say that $S_{1}$ and $S_{2}$ gives 
a \emph{{\bf splitting of $n$ by $\mu_{b}$ over $\FF_{q}$}}. This gives us another way of describing 
duadic codes.  Note that for a fixed pair of duadic codes over $\FF_{q}$ of length $n$, we can use the 
same multiplier for the splitting in Definition \ref{duadicdef} and the splitting of $n$  in Theorem 
\ref{Thm 6.1.5}.

\begin{thm}
{\rm (\cite{fundamentals}.)}
Duadic codes of length $n$ over $\FF_{q}$ exist if and only if $q$ is a square mod~$n$.
\label{Thm 6.3.2}
\end{thm}

\subsection{Hermitian Orthogonality of Duadic Codes over $\FF_{q^{2}}$}

From this point onwards, we consider codes over the Galois field $\FF_{q^{2}}$, where $q$ is a power of 
some prime $p$. Again we assume that $n$ is an odd positive integer and gcd$(n,q)=1$. Thus duadic codes of 
length $n$ over $\FF_{q^{2}}$ always exist by Theorem \ref{Thm 6.3.2}. The following theorem is the 
Hermitian analogue of Theorem 6.4.1 of \cite{fundamentals}, where the Euclidean self-orthogonality of 
duadic codes over $\FF_{q}$ are considered.

\begin{prop}
Let $C$ be any $\lbrack n, \frac{n-1}{2} \rbrack$ cyclic code of length $n$ over $\FF_{q^{2}}$. Then $C$ 
is Hermitian self-orthogonal if and only if $C$ is an even-like duadic code whose splitting is given by 
$\mu_{-q}$.
\label{Thm 6.4.1.1}
\end{prop}

\pf $(\Leftarrow)$ Suppose $C = C_{1}$ is an even-like duadic code whose splitting is given by $\mu_{-q}$. 
Let $e(x)$ be the generating idempotent for $C$. By Theorem 6.1.3 (vi) of \cite{fundamentals}, $C = C_{1} 
\subseteq D_{1} = \langle 1 - e(x)\mu_{-q} \rangle$.
By Proposition \ref{Thm 4.4.9.1}, the generating idempotent for $C^{\bot_{H}}$ is also $1-e(x)\mu_{-q}$.
 Thus $D_{1} = C^{\bot_{H}}$ and so $C$ is Hermitian self-orthogonal.
\vskip .1in
$(\Rightarrow)$ Let $C = C_{1}$ be a Hermitian self-orthogonal cyclic code. Let $e_{1}(x)$ be the 
generating idempotent for $C_{1}$ and $T_{1}$ its defining set. Since $C_{1}$ is Hermitian self-orthogonal 
and $\overline{j}(x)$ is not orthogonal to itself, $\overline{j}(x) \not\in C_{1}$. Hence by Lemma 6.1.2 
(iii) of \cite{fundamentals}, $C_{1}$ is even-like.
Let $e_{2}(x) = e_{1}(x)\mu_{-q}$ and let $C_{2} = \langle e_{2}(x) \rangle$. By Theorem 4.3.13 of 
\cite{fundamentals}, $C_{2} = C_{1}\mu_{-q}$.
\vskip 0in
Let $(a_{0},a_{1},\ldots,a_{n-1}) \in C_{1}$. 
Since $C_{1}$ is even-like, it follows that $\sum_{i=0}^{n-1}a_{i} = 0$. Thus
$(1,1,\ldots,1)\cdot \overline{(a_{0},a_{1},\ldots,a_{n-1})} 
 =  (1,1,\ldots,1) \cdot (a_{0}^{q},a_{1}^{q},\ldots,a_{n-1}^{q})
 =  \sum_{i=0}^{n-1}a_{i}^{q}
 =  \left( \sum_{i=0}^{n-1}a_{i} \right)^{q}~=~0$
which implies that $\overline{j}(x) \in C_{1}^{\bot_{H}}$.
Since $C_{1}^{\bot_{H}}$ has dimension $\frac{n+1}{2}$ and $C_{1} \subseteq C_{1}^{\bot_{H}}$, we have 
$C_{1}^{\bot_{H}} = C_{1} + \langle \overline{j}(x) \rangle$. 
Using Theorem 4.3.7 of \cite{fundamentals} and Lemma 6.1.2 (i) of \cite{fundamentals}, $C_{1}^{\bot_{H}}$ 
has generating idempotent $e_{1}(x) + \overline{j}(x)$.
By Proposition \ref{Thm 4.4.9.1}, the generating idempotent for $C_{1}^{\bot_{H}}$ is $1 - 
e_{1}(x)\mu_{-q}$. By the uniqueness of the idempotent generator, we must have
 $ 1 - e_{1}(x)\mu_{-q}  =  e_{1}(x) + \overline{j}(x)$
which implies  $1 - \overline{j}(x) = e_{1}(x) + e_{1}(x)\mu_{-q} = e_{1}(x) + e_{2}(x)$.
 Clearly $e_{1}(x) = e_{2}(x)(\mu_{-q})^{-1} = e_{2}(x)(\mu_{-q})$.
 Therefore $C_{1}$ and $C_{2}$ form a pair of even-like codes whose splitting is given by $\mu_{-q}$. 
\hfill \Q{\Box}

\vskip .1in
\begin{lem}
Let $C$ be a cyclic code. Then $(C\mu_{a})^{\bot_{H}} = C^{\bot_{H}}\mu_{a}$.
\label{LEMMA}
\end{lem}

\pf  Use Proposition \ref{Thm 4.4.9.1} above and Theorem 4.3.13 of \cite{fundamentals} to show that 
$(C\mu_{a})^{\bot_{H}}$ and $C^{\bot_{H}}\mu_{a}$ have the same idempotent generator.
\hfill \Q{\Box}

\vskip .1in
\begin{prop}
Suppose that $C_{1}$ and $C_{2}$ are a pair of even-like duadic codes over $\FF_{q^{2}}$, having $D_{1}$ 
and $D_{2}$ as their associated odd-like duadic codes. Then the following are equivalent.
\begin{enumerate}
\item $C_{1}^{\bot_{H}} = D_{1}$
\item $C_{2}^{\bot_{H}} = D_{2}$
\item $C_{1}\mu_{-q} = C_{2}$
\item $C_{2}\mu_{-q} = C_{1}$
\end{enumerate}
\label{Thm 6.4.2.1}
\end{prop}
 
\pf  From the definition of duadic codes and Theorem  6.1.3 (vii) of \cite{fundamentals}, we obtain 
$C_{1}\mu_{a} = C_{2}$, $C_{2}\mu_{a} = C_{1}$, $D_{1}\mu_{a} = D_{2}$ and $D_{2}\mu_{a} = D_{1}$ for some 
$a$. Hence by Lemma \ref{LEMMA},
\vskip 0in
\noindent if part 1 holds, then
\begin{equation*}
 C_{2}^{\bot_{H}} = (C_{1}\mu_{a})^{\bot_{H}} = C_{1}^{\bot_{H}}\mu_{a} = D_{1}\mu_{a} = D_{2}
\end{equation*}
and if part 2 holds,  then
\begin{equation*}
 C_{1}^{\bot_{H}} = (C_{2}\mu_{a})^{\bot_{H}} = C_{2}^{\bot_{H}}\mu_{a} = D_{2}\mu_{a} = D_{1}.
\end{equation*}
Hence parts 1 and 2 are equivalent.

Part 3 is equivalent to part 4 since $(\mu_{-q})^{-1} = \mu_{-q}$.

If part 1 holds, then by  Theorem 6.1.3 (vi) of \cite{fundamentals}, $C_{1}$ is Hermitian self-orthogonal. 
Hence by Proposition \ref{Thm 6.4.1.1}, part 3 holds.

If part 3 holds, then $\mu_{-q}$ gives a splitting for $C_{1}$ and $C_{2}$. Let $e_{i}(x)$ be the 
generating idempotent for $C_{i}$. By Theorem 4.3.13 of \cite{fundamentals}, $e_{1}(x)\mu_{-q} = 
e_{2}(x)$. Hence by Proposition \ref{Thm 4.4.9.1}, the generating idempotent for $C_{1}^{\bot_{H}}$ is $1 
- e_{1}(x)\mu_{-q} = 1 - e_{2}(x)$. Thus, part 1 holds, completing the proof. \hfill \Q{\Box}

\vskip .1in

\begin{prop}
Suppose that $C_{1}$ and $C_{2}$ are a pair of even-like duadic codes over $\FF_{q^{2}}$, having $D_{1}$ 
and $D_{2}$ as their associated odd-like duadic codes. Then the following are equivalent.
\begin{enumerate}
\item $C_{1}^{\bot_{H}} = D_{2}$
\item $C_{2}^{\bot_{H}} = D_{1}$
\item $C_{1}\mu_{-q} = C_{1}$
\item $C_{2}\mu_{-q} = C_{2}$
\end{enumerate}
\label{Thm 6.4.3.1}
\end{prop}

\pf From the definition of duadic codes and  Theorem  6.1.3 (vii) of \cite{fundamentals}, we obtain 
$C_{1}\mu_{a} = C_{2}$, $C_{2}\mu_{a} = C_{1}$, $D_{1}\mu_{a} = D_{2}$ and $D_{2}\mu_{a} = D_{1}$ for some 
$a$. Hence, by Lemma \ref{LEMMA},
\vskip 0in
\noindent if part 1 holds, then
\begin{equation*}
 C_{2}^{\bot_{H}} = (C_{1}\mu_{a})^{\bot_{H}} = C_{1}^{\bot_{H}}\mu_{a} = D_{2}\mu_{a} = D_{1}
\end{equation*}
and if part 2 holds,  then
\begin{equation*}
 C_{1}^{\bot_{H}} = (C_{2}\mu_{a})^{\bot_{H}} = C_{2}^{\bot_{H}}\mu_{a} = D_{1}\mu_{a} = D_{2}.
\end{equation*}
Hence parts 1 and 2 are equivalent.

Let $e_{i}(x)$ be the generating idempotent for $C_{i}$. By Proposition \ref{Thm 4.4.9.1}, 
$C_{1}^{\bot_{H}}$ has generating idempotent $1 - e_{1}(x)\mu_{-q}$. Thus
$C_{1}^{\bot_{H}} = D_{2}$ if and only if  $1 - e_{1}(x)\mu_{-q} = 1 - e_{1}(x)$ if and only if  
$e_{1}(x)\mu_{-q} = e_{1}(x)$ if and only if $C_{1}\mu_{-q} = C_{1}$ by Theorem 4.3.13 of 
\cite{fundamentals}. Hence parts 1 and 3 are equivalent. It can be shown by an analogous argument that 
parts 2 and 4 are equivalent. \hfill \Q{\Box}

\vskip .1in

\subsection{Extensions of Odd-like Duadic Codes}

Odd-like duadic codes have parameters $\lbrack n , \frac{n+1}{2} \rbrack$. Hence it is interesting to 
consider extending such codes because such extensions could possibly be Hermitian self-dual codes. The 
goal of this section is to give a way of extending odd-like duadic codes and to give conditions under 
which these extensions are Hermitian self-dual. We also prove that any cyclic code whose extended code is 
Hermitian self-dual must be an odd-like duadic code.

Let $D$ be an odd-like duadic code.  The code $D$ can be obtained from its even-like subcode $C$ by adding 
$\overline{j}(x)$ to a basis of $C$ (Theorem 6.1.3 (ix), \cite{fundamentals} ). Hence it is natural to 
define an extension for which the all-one vector $\mathbf{1}$ is Hermitian orthogonal to itself.

\noindent In $\FF_{q^{2}}$ consider the equation 
\begin{equation}
1 + \gamma^{q+1}n = 0.
\label{gamma1}
\end{equation}

\noindent Since $q$ is a power of a prime $p$ and  $n \in \FF_{p} \subseteq \FF_{q}$, we have $n^{q} = n$, 
or $n^{q+1} = n^{2}$ in $\FF_{q^{2}}$. So
$1 + \gamma^{q+1}n = 0 \iff  n + \gamma^{q+1}n^{2} = 0 \iff  n + \gamma^{q+1}n^{q+1} = 0 \iff  n + (\gamma 
n)^{q+1} = 0$. Thus Equation (\ref{gamma1}) is equivalent to
\begin{equation}
n + \gamma^{q+1} = 0.
\label{gamma2}
\end{equation}
Note that $\lbrace a^{q+1} \mid a \in \FF_{q^{2}} \rbrace = \FF_{q}$.
Thus (\ref{gamma2}) will always have a solution in $\FF_{q^{2}}$, which implies that (\ref{gamma1}) is 
solvable in $\FF_{q^{2}}$.

We are now ready to describe the extension.
Let $\gamma$ be a solution to (\ref{gamma1}).
Let $\mathbf{c} = ( c_{0},c_{1},\ldots,c_{n-1} ) \in D$.
Define the extended codeword $\widetilde{\mathbf{c}} = ( c_{0},c_{1},\ldots,c_{n-1},c_{\infty} )$, 
where
\begin{equation*}
c_{\infty} = -\gamma\sum_{i=0}^{n-1}c_{i}.
\end{equation*}
Let $\widetilde{D} = \lbrace \widetilde{\mathbf{c}} \mid \mathbf{c} \in D \rbrace $ be the extended code 
of $D$.

\begin{prop}
Let $D_{1}$ and  $D_{2}$ be a pair of odd-like duadic codes of length $n$ over $\FF_{q^{2}}$. The 
following hold:
\begin{enumerate}
\item If $\mu_{-q}$ gives the splitting for $D_{1}$ and  $D_{2}$, then $\widetilde{D_{1}}$ and  
$\widetilde{D_{2}}$ are Hermitian self-dual.
\item If $D_{1}\mu_{-q} = D_{1}$, then $\widetilde{D_{1}}$ and  $\widetilde{D_{2}}$ are Hermitian duals of 
each other.
\end{enumerate}
\label{Thm 6.4.14.1}
\end{prop}

\pf Let $C_{1}$ and  $C_{2}$ be the even-like duadic codes associated to $D_{1}$ and  $D_{2}$.

\noindent Note that
\begin{eqnarray*}
\widetilde{\overline{j}(x)} \overline{\widetilde{\overline{j}(x)}} & = & \left( \frac{1}{n}, 
\frac{1}{n},\ldots , \frac{1}{n}, -\gamma \right) \cdot \overline{\left( \frac{1}{n}, \frac{1}{n},\ldots , 
\frac{1}{n}, -\gamma \right)} \\
& = & \left( \frac{1}{n}, \frac{1}{n},\ldots , \frac{1}{n}, -\gamma \right) \cdot \left( \frac{1}{n}, 
\frac{1}{n},\ldots , \frac{1}{n}, (-\gamma)^{q} \right) \\
& = & \frac{1}{n} + \gamma^{q+1} \\
& = & \frac{1}{n}(1 + \gamma^{q+1}n) \\
& = & 0,
\end{eqnarray*}
by our choice of $\gamma$. This shows that $\widetilde{\overline{j}(x)}$ is Hermitian orthogonal to 
itself. Since $C_{i}$ is even-like, $\widetilde{C_{i}}$ is obtained by adding a zero coordinate to $C_{i}$ 
and so $\widetilde{\overline{j}(x)}$ is also orthogonal to $\widetilde{C_{i}}$.

We first prove part 1. Proposition \ref{Thm 6.4.1.1} ensures that $C_{1}$ is Hermitian self-orthogonal, 
and so $\widetilde{C_{1}}$ is Hermitian self-orthogonal.
Since $\widetilde{\overline{j}(x)}$ is orthogonal to $\widetilde{C_{1}}$, the code spanned by $\langle 
\widetilde{C_{1}}, \widetilde{\overline{j}(x)} \rangle$ is Hermitian self-orthogonal.
However, by  Theorem 6.1.3 (ix) of \cite{fundamentals}, $D_{1} = \langle C_{1}, \overline{j}(x) \rangle$. 
Clearly $\widetilde{D_{1}} = \langle \widetilde{C_{1}}, \widetilde{\overline{j}(x)} \rangle$. Thus 
$\widetilde{D_{1}}$ is Hermitian self-orthogonal. Since the dimension of $\widetilde{D_{1}}$ is 
$\frac{n+1}{2}$, $\widetilde{D_{1}}$ is Hermitian self dual.
Analogous arguments will prove that $\widetilde{D_{2}}$ is Hermitian self-dual.
\vskip 0in
We now prove part 2. Suppose $D_{1}\mu_{-q} = D_{1}$. It follows that $C_{1}\mu_{-q}=C_{1}$. 
By Proposition \ref{Thm 6.4.3.1}, $C_{2}^{\bot_{H}} = D_{1}$ and so by Theorem  6.1.3  (vi) of 
\cite{fundamentals}, $ C_{1} \subseteq C_{2}^{\bot_{H}}$.
Therefore $\widetilde{C_{1}}$ and $\widetilde{C_{2}}$ are orthogonal to each other and consequently the 
codes spanned by $\langle \widetilde{C_{1}}, \widetilde{\overline{j}(x)} \rangle$ and $\langle 
\widetilde{C_{2}}, \widetilde{\overline{j}(x)} \rangle$ are orthogonal.
By  Theorem 6.1.3 (v) \& (vi) of  \cite{fundamentals}, these codes must be $\widetilde{D_{1}}$ and  
$\widetilde{D_{2}}$ of dimension $\frac{n+1}{2}$. Therefore $\widetilde{D_{1}}$ and  $\widetilde{D_{2}}$ 
are duals of each other. \hfill \Q{\Box}

\begin{cor}
Let $C$ be a cyclic code over $\FF_{q^{2}}$. The extended code $\widetilde{C}$ is Hermitian self-dual if 
and only if $C$ is an odd-like duadic code whose splitting is given by $\mu_{-q}$.
\label{Cor 6.4.14.1}
\end{cor}

\pf  

\noindent $(\Leftarrow)$ This follows directly from the preceding proposition.

\noindent $(\Rightarrow)$ Since the extended code $\widetilde{C}$ has length $n+1$, the dimension of $C$ 
is $\frac{n+1}{2}$ and therefore $C$ cannot be Hermitian self-orthogonal.
The assumption that $\widetilde{C}$ is self-dual implies that the even-like subcode of $C$ is necessarily 
Hermitian self-orthogonal. Since $C$ is not Hermitian self-orthogonal, $C$ cannot be even-like. Let 
$C_{e}$ be the even-like subcode of $C$. The code $C_{e}$ is an $[n, \frac{n-1}{2}]$ Hermitian 
self-orthogonal cyclic code and so by Proposition \ref{Thm 6.4.1.1}, $C_{e}$ is an even-like duadic code 
with splitting by $\mu_{-q}$. Thus $C$ is an odd-like duadic code with a splitting by $\mu_{-q}$.
\hfill \Q{\Box}

\vskip .1in

\section{Lengths with Splittings By $\mu_{-q}$}

All throughout this section, we let $q$ be a power of a prime $p$ and we assume that $n$ is an odd integer 
with $gcd(n,q) = 1$. Define $ord_{r}(q)$ to be the smallest positive integer $t$ such that $q^{t} \equiv 1 
\;(\mbox{mod }~r)$. 
In view of Proposition \ref{Thm 6.4.1.1} and Corollary \ref{Cor 6.4.14.1}, it is natural to ask under what 
conditions do we get a splitting of $n$ by $\mu_{-q}$. We note that the study of the feasibility of an 
integer in \cite{SD over GF4} becomes a special case of this with $q=2$.
\vskip .1in
\noindent The main result of this section is the following theorem.

\begin{thm}
The permutation map $\mu_{-q}$ gives a splitting of $n$ if and only if $ord_{r}(q) \not\equiv 2  \; 
(\mbox{mod } 4)$ for every prime $r$ dividing $n$.
\label{splitting}
\end{thm}

\noindent Our proof of this theorem will be based on several lemmas. Lemma \ref{prop 3} is a well-known 
fact from elementary number theory, see e.g. Proposition 3 in \cite{moree}, and we leave its proof as an 
exercise to the reader.

\begin{lem}
Let $r$ be a prime distinct from $p$. Then $r$ divides $q^{k} + 1$ for some positive integer $k$ if and 
only if $ord_{r}(q)$ is even.
\label{prop 3}
\end{lem}

\begin{lem}
Let $r$ be a prime distinct from $p$. Then $r$ divides $q^{2i-1}+1$ for some integer $i \geq 1 $ if and 
only if $ord_{r}(q) \equiv 2 \;(\mbox{mod } 4)$. 
\label{Lem a}
\end{lem}

\pf By Lemma \ref{prop 3}, $r$ divides $q^{k}+1$ for some positive integer $k$ if and only if  
$ord_{r}(q)$ is even.
If $ord_{r}(q)$ is even, then
\begin{equation*}
r \mid q^{k}+1 \Longleftrightarrow q^{k} \equiv -1 \;(\mbox{mod } r) \Longleftrightarrow  
k \equiv \frac{ord_{r}(q)}{2} \; (\mbox{mod } ord_{r}(q)). 
\end{equation*}
Thus $r$ divides $q^{2i-1}+1$ if and only if $ord_{r}(q)$ is even and 
\begin{equation}
2i-1 \equiv \frac{ord_{r}(q)}{2} \; (\mbox{mod } ord_{r}(q)).
\label{sol}
\end{equation}
But (\ref{sol}) has a solution $i$ if and only if $ord_{r}(q) \equiv 2  \;(\mbox{mod } 4)$. \hfill 
\Q{\Box}

\begin{prop}
Assume $gcd(n,q)=1$. Then
$gcd(n,q^{2i-1}+1) = 1$ for every integer $i \geq 1$ if and only if $ord_{r}(q) \not\equiv 2 \;(\mbox{mod 
} 4)$ for every prime $r$ dividing $n$.
\label{Prop 1}
\end{prop}

\pf  Write $n = r_{1}^{e_{1}} r_{2}^{e_{2}} \cdots r_{s}^{e_{s}}$. Then, using Lemma \ref{Lem a},
$ord_{r_{j}}(q) \not\equiv 2 \;(\mbox{mod } 4) \mbox{ for all } j = 1,\ldots,s$ if and only if for all $j 
= 1,\ldots,s, \; \; r_{j} \mbox{ does not divide } q^{2i-1}+1 \mbox{ for every } i \geq 1$ if and only if 
for all $j = 1,2,\ldots,s, \; \; gcd(r_{j}, q^{2i-1}+1) = 1 \mbox{ for every } i \geq 1$ if and only if 
$gcd(n,q^{2i-1}+1) = 1  \mbox{ for every } i \geq 1$.

\begin{prop}
Let $t$ be an integer such that  $t \not\equiv (q^{2})^{j} \;(\mbox{mod } n)$ and $t^{2} \equiv 
(q^{2})^{j} \;(\mbox{mod } n)$ for some non-negative integer $j$. Suppose $gcd(t,n)=1$. Then $\mu_{t}$ 
gives a splitting of $n$ if and only if $gcd(n,q^{2i}-t) = 1$ for every integer $i \geq 1$.
\label{Prop 2}
\end{prop}

\pf  Clearly by the assumptions on $t$, $(\mu_{t})^{2}(C_{s}) = C_{s}$ for every $q^{2}$-cyclotomic coset 
$C_{s}$.
Thus $\mu_{t}$ gives a splitting of $n$ if and only if it does not fix any $q^{2}$-cyclotomic coset.
Let $C_{a}$ be a $q^{2}$-cyclotomic coset. Then $\mu_{t}$ fixes $C_{a}$ if and only if $ta \equiv 
(q^{2})^{i}a \;(\mbox{mod } n)$ for some positive integer $i$. Thus $\mu_{t}$ gives a splitting of $n$ if 
and only if $ta \not\equiv (q^{2})^{i}a \;(\mbox{mod } n)$ for every $i \geq 1$ if and only if 
$gcd(n,q^{2i}-t) = 1 \mbox{ for every } i \geq 1$. \hfill \Q{\Box}

\vskip .1in
\noindent Theorem 9 of \cite{qcodes} is a special case of Proposition \ref{Prop 2} with $q=2$.

\begin{cor}
The permutation map $\mu_{-q}$ gives a splitting of $n$ if and only if $gcd(n,q^{2i-1}+1) = 1$ for every 
integer $i \geq 1$.
\label{Cor 2}
\end{cor}

\pf  This follows immediately from Proposition \ref{Prop 2} since $gcd(n,q) = 1$ by assumption. \hfill 
\Q{\Box}

\vskip .1in
\noindent We are now ready to prove the main theorem of this section.
\vskip .1in
\noindent {\bf Proof of Theorem \ref{splitting}:}
\nopagebreak

By Corollary \ref{Cor 2}, the permutation map $\mu_{-q}$ gives a splitting of $n$ if and only if 
$gcd(n,q^{2i-1}+1) = 1$ for every integer $i \geq 1$.
By Proposition \ref{Prop 1}, $gcd(n,q^{2i-1}+1) = 1$ for every integer $i \geq 1$ if and only if 
$ord_{r}(q) \not\equiv 2 \;(\mbox{mod } 4)$ for every prime $r$ dividing $n$.  \hfill \Q{\Box}

\vskip .1in
\noindent We remark that Theorem \ref{splitting} says that $\mu_{-q}$ gives a splitting of $n$ if and only 
if for all prime $r$ dividing $n$,  either $ord_{r}(q)$ is odd or $ord_{r}(q)$ is doubly even. However, it 
is easy to show that $ord_{r}(q)$ is doubly even if and only if $ord_{r}(q^{2})$ is even.
Thus we can restate Theorem \ref{splitting} as:

\begin{thm}
The permutation map $\mu_{-q}$ gives a splitting of $n$ if and only if for every prime $r$ dividing $n$,  
either $ord_{r}(q)$ is odd or $ord_{r}(q^{2})$ is even.
\label{Thm 3.17}
\end{thm}

Lastly, we arrive at the following result which gives sufficient and necessary conditions for the 
existence of  Hermitian self-dual extended cyclic code. We note that the same result was obtained in 
\cite{conchita} for the more general case of group codes.

\begin{thm}
Cyclic codes of length $n$ over $\FF_{q^{2}}$ whose extended code is Hermitian self-dual exist if and only 
if for every prime $r$ dividing $n$,  either $ord_{r}(q)$ is odd or $ord_{r}(q^{2})$ is even.
\label{concon}
\end{thm}

\pf This follows directly from Corollary \ref{Cor 6.4.14.1} and Theorem \ref{Thm 3.17}.
\hfill \Q{\Box}

\vskip .1in

The table below enumerates all the splittings (up to symmetry between $S_{1}$ and $S_{2}$) of $n$ by 
$\mu_{-q}$ over $\FF_{q^{2}}$ for $5 \leq n \leq 45$ and $q=3$, $q=4$ and $q=5$ by listing all the 
possible sets for the $S_{1}$ in Theorem \ref{Thm 6.1.5}. The $C_{i}$'s are $q^{2}$-cyclotomic cosets 
modulo $n$. We omit those $n$ for which no such splitting exists for all values of $q$.

\begin{table}
\scriptsize
\begin{center}
\begin{tabular}{|c|l|l|l|}  \hline
   \bf{$n$} &       \multicolumn{3}{c|}{\bf{$S_{1}$}}  \\
		                   \hline
       &     \hspace{.5in}   \bf{$q=3$}   &   \hspace{1.2in}  \bf{$q=4$ } & \hspace{.65in}  \bf{$q=5$} \\
		   					                \hline
    5  &       $C_{1}^{\; \; \clubsuit}$ & \hspace{.65in}  ------ & \hspace{.65in}  ------  \\ \hline
    7  &      \hspace{.65in}  ------ & $C_{1}^{\; \; \clubsuit}$  & \hspace{.65in}  ------ \\ \hline
	9  &      \hspace{.65in}  ------ & $C_{1} \cup C_{3}$, $C_{1} \cup C_{6}$  & \hspace{.65in}  ------ \\ 
\hline
    11 &      $C_{1}^{\; \; \clubsuit}$ & $C_{1}^{\; \; \clubsuit}$ & $C_{1}^{\; \; \clubsuit}$ \\ \hline
    13 &      $C_{1} \cup C_{2}$, $C_{1} \cup C_{7}$ & \hspace{.65in}  ------ & $C_{1} \cup C_{2} \cup 
C_{4}$, $C_{1} \cup C_{2} \cup C_{6}$,  \\
       &  &  & $C_{1} \cup C_{3} \cup C_{4}^{\; \; \; \clubsuit}$, $C_{1} \cup C_{3} \cup C_{6}$   \\ 
\hline
    17 &      $C_{1}^{\; \; \clubsuit}$ & $C_{1} \cup C_{2} \cup C_{3} \cup C_{6}$, $C_{1} \cup C_{2} \cup  		
         C_{3} \cup C_{7}$,  & $C_{1}^{\; \; \clubsuit}$ \\ 
       &  &  $C_{1} \cup C_{2} \cup C_{5} \cup C_{6}$, $C_{1} \cup C_{2} \cup C_{5} \cup C_{7}$,  & \\ 
	   &  &  $C_{1} \cup C_{3} \cup C_{6} \cup C_{8}$, $C_{1} \cup C_{3} \cup C_{7} \cup C_{8}$,  & \\
	   &  &  $C_{1} \cup C_{5} \cup C_{6} \cup C_{8}$, $C_{1} \cup C_{5} \cup C_{7} \cup C_{8}$, & 
				\\ 	\hline
	19 & \hspace{.65in}  ------  &   $C_{1}^{\; \; \clubsuit}$  & $C_{1}^{\; \; \clubsuit}$ \\ \hline
	21 & \hspace{.65in}  ------ & $C_{1} \cup C_{2} \cup C_{3} \cup C_{7}$, $C_{1} \cup C_{2} \cup C_{3} 
				\cup C_{14}$, & \hspace{.65in}  ------ \\
		&  & $C_{1} \cup C_{2} \cup C_{7} \cup C_{9}$, $C_{1} \cup C_{2} \cup C_{9} \cup C_{14}$, & \\
		&  & $C_{1} \cup C_{3} \cup C_{7} \cup C_{10}$, $C_{1} \cup C_{3} \cup C_{10} \cup C_{14}$, & \\
		&  & $C_{1} \cup C_{7} \cup C_{9} \cup C_{10}$, $C_{1} \cup C_{9} \cup C_{10} \cup C_{14}$ & \\
	\hline
    23 &  $C_{1}^{\; \; \clubsuit}$ & $C_{1}^{\; \; \clubsuit}$ & \hspace{.65in}  ------ \\ \hline
    25 &      $C_{1} \cup C_{5}$, $C_{1} \cup C_{10}$ & \hspace{.65in}  ------ & \hspace{.65in}  ------ \\ 
\hline
	27 & \hspace{.65in}  ------ & $C_{1} \cup C_{3} \cup C_{9}$, $C_{1} \cup C _{3} \cup C_{18}$, $C_{1} 
\cup C_{6} \cup C_{9}$, $C_{1} \cup C _{6} \cup C_{18}$ & 
			\hspace{.65in}  ------ \\
	  \hline
    29 &  $C_{1}^{\; \; \clubsuit}$ & \hspace{.65in}  ------ & \hspace{.65in}  ------ \\ \hline
	31 &  \hspace{.65in} ------ & $C_{1} \cup C_{3} \cup C_{5}$, $C_{1} \cup C_{3} \cup C_{11}$, & $C_{1} 
\cup C_{2} \cup C_{3} \cup C_{4} \cup C_{8}$,  \\		    
       &  & $C_{1} \cup C_{5} \cup C_{7}^{\; \; \clubsuit}$, $C_{1} \cup C_{7} \cup C_{11}$,  & $C_{1} 
\cup C_{2} \cup C_{3} \cup C_{4} \cup C_{17}$, \\
	   &  &  & $C_{1} \cup C_{2} \cup C_{3} \cup C_{8} \cup C_{11}$, \\
	   &  &  & $C_{1} \cup C_{2} \cup C_{3} \cup C_{11} \cup C_{17}$, \\
	   &  &  & $C_{1} \cup C_{2} \cup C_{4} \cup C_{8} \cup C_{16}^{\; \; \clubsuit}$, \\
	   &  &  & $C_{1} \cup C_{2} \cup C_{4} \cup C_{16} \cup C_{17}$, \\
	   &  &  & $C_{1} \cup C_{2} \cup C_{8} \cup C_{11} \cup C_{16}$, \\
	   &  &  & $C_{1} \cup C_{2} \cup C_{11} \cup C_{16} \cup C_{17}$, \\
	   &  &  & $C_{1} \cup C_{3} \cup C_{4} \cup C_{8} \cup C_{12}$, \\
	   &  &  & $C_{1} \cup C_{3} \cup C_{4} \cup C_{12} \cup C_{17}$, \\
	   &  &  & $C_{1} \cup C_{3} \cup C_{8} \cup C_{11} \cup C_{12}$, \\
	   &  &  & $C_{1} \cup C_{3} \cup C_{11} \cup C_{12} \cup C_{17}$, \\
	   &  &  & $C_{1} \cup C_{4} \cup C_{8} \cup C_{12} \cup C_{16}$, \\
	   &  &  & $C_{1} \cup C_{4} \cup C_{12} \cup C_{16} \cup C_{17}$, \\
	   &  &  & $C_{1} \cup C_{8} \cup C_{11} \cup C_{12} \cup C_{16}$, \\
	   &  &  & $C_{1} \cup C_{11} \cup C_{12} \cup C_{16} \cup C_{17}$, \\
		 \hline
	33 & \hspace{.65in} ------ & $C_{1} \cup C_{3} \cup C_{5} \cup C_{11}$, $C_{1} \cup C_{3} \cup C_{5} 
			\cup C_{22}$, & \hspace{.65in} ------ \\
       &  &  $C_{1} \cup C_{3} \cup C_{7} \cup C_{11}$, $C_{1} \cup C_{3} \cup C_{7} \cup C_{22}$, & \\
	   &  &  $C_{1} \cup C_{5} \cup C_{6} \cup C_{11}$, $C_{1} \cup C_{5} \cup C_{6} \cup C_{22}$, & \\
	   &  &  $C_{1} \cup C_{6} \cup C_{7} \cup C_{11}$, $C_{1} \cup C_{6} \cup C_{7} \cup C_{22}$, & \\
	   			\hline
	37 &  \hspace{.65in} ------ & \hspace{.65in} ------ & $C_{1}^{\; \; \; \clubsuit}$ \\ \hline
	41 &  $C_{1} \cup C_{2}\cup C_{4}\cup C_{7}\cup C_{8}$, & \hspace{.65in} ------ & $C_{1} \cup C_{3}$, 
$C_{1} \cup C_{6}$ \\
	   &  $C_{1} \cup C_{2}\cup C_{4}\cup C_{7}\cup C_{11}$, & & \\ 
       &  $C_{1} \cup C_{2}\cup C_{4}\cup C_{8}\cup C_{16}^{\; \; \; \clubsuit}$, & & \\
	   &  $C_{1} \cup C_{2}\cup C_{4}\cup C_{11}\cup C_{16}$, & & \\ 
	   &  $C_{1} \cup C_{2}\cup C_{7}\cup C_{8}\cup C_{12}$, & & \\
	   &  $C_{1} \cup C_{2}\cup C_{7}\cup C_{11}\cup C_{12}$, & & \\ 
	   &  $C_{1} \cup C_{2}\cup C_{8}\cup C_{12}\cup C_{16}$, & & \\
	   &  $C_{1} \cup C_{2}\cup C_{11}\cup C_{12}\cup C_{16}$, & & \\ 
	   &  $C_{1} \cup C_{4}\cup C_{6}\cup C_{7}\cup C_{8}$, & & \\
	   &  $C_{1} \cup C_{4}\cup C_{6}\cup C_{7}\cup C_{11}$, & & \\ 
	   &  $C_{1} \cup C_{4}\cup C_{6}\cup C_{8}\cup C_{16}$, & & \\
	   &  $C_{1} \cup C_{4}\cup C_{6}\cup C_{11}\cup C_{16}$, & & \\ 
	   &  $C_{1} \cup C_{6}\cup C_{7}\cup C_{8}\cup C_{12}$, & & \\
	   &  $C_{1} \cup C_{6}\cup C_{7}\cup C_{11}\cup C_{12}$, & & \\ 
	   &  $C_{1} \cup C_{6}\cup C_{8}\cup C_{12}\cup C_{16}$, & & \\
	   &  $C_{1} \cup C_{6}\cup C_{11}\cup C_{12}\cup C_{16}$, & & \\ 
												   \hline
	43 & \hspace{.65in} ------ & $C_{1} \cup C_{3} \cup C_{7}$, $C_{1} \cup C_{3} \cup C_{9}$, $C_{1} \cup 
C_{6} \cup C_{7}$, $C_{1} \cup C_{6} \cup C_{9}^{\; \; \clubsuit}$& 
				\hspace{.65in} ------ \\
           \hline

\end{tabular}
\caption{Splittings of $n$ by $\mu_{-q}$ \footnotesize{($\clubsuit$ denotes splittings of Quadratic 
Residue codes)}}				   
\end{center}
\end{table}


\appendix

\section{Quantitative Aspects}

\rm

\subsection{Counting integers that are split by $\mu_{-q}$}

Theorem \ref{splitting} raises the question of counting the number of
integers $n\le x$ such that $\mu_{-q}$ gives a splitting
of $n$.
In other words, we are interested in counting those integers $n$ such
that $n$ is coprime with the sequence $S(q):=\{q^{2i-1}+1\}_{i=1}^{\infty}$. 
We let $A_q(x)$ denote the
associated counting function. We are interested in sharp estimates
for $A_q(x)$ as $x$ gets large. We use the shorthand GRH to denote
the Generalized Riemann Hypothesis. The best we can do in this
respect is stated in the following theorem:

\begin{thm} 
\label{main}
Let $q=p^{t}$ be a prime power. Put $\lambda=\nu_{2}(t)$.
\begin{enumerate}
\item For some positive constant $c_{q}$ we have
\begin{equation*}
A_{q}(x)=c_{q}{x \over \log^{\delta(q)}x}
+O_{q}\left({x(\log \log x)^5\over \log^{1+\delta(q)}x}\right),
\end{equation*}
where the implicit constant depends at most on $q$.
\item Let $\epsilon>0$ and $v\ge 1$ be
arbitrary. Assuming GRH we have that
\begin{equation*}
A_{q}(x)=
\sum_{0\le j <v}{b_{j}x\over \log^{\delta(q)+j}x}+O_{\epsilon,q}\left({x\over
\log^{\delta(q)+v-\epsilon}x}\right),
\end{equation*}
where the implied constant depends at most on $\epsilon$ and $q$, and
$b_0(=c_{q}),\ldots,b_{v}$ are constants that depend at most on $q$.
\end{enumerate}

The constant $\delta(q)$ is the natural density of primes $r$ such
that ord$_{r}(q)\equiv 2 \; (\mbox{mod } 4)$ and is given as follows:
\begin{equation*}
\delta(p^{t})=\left\{ 
\begin{array}{ll}
7/24                & \mbox{ if }  p=2 \mbox{ and }  \lambda=0 ; \\
1/3                 & \mbox{ if } p=2 \mbox{ and } \lambda=1 ; \\
2^{-\lambda -1}/3   & \mbox{ if } p=2 \mbox{ and } \lambda\ge 2 ;\\
2^{-\lambda}/3      & \mbox{ if } p\ne 2.
\end{array}
\right.
\end{equation*}
\end{thm}

\noindent Our proof of Theorem \ref{main} rests on various lemmas.
Let $\chi_{q}(n)$ be the characteristic function of the integers
$n$ that are coprime with the sequence $S(q)$, i.e.
\begin{equation*}
\chi_{q}(n)=\left\{ 
\begin{array}{ll}
1 & \mbox{ if } (n,S(q))=1 ; \\
0 & \mbox{ otherwise}. 
\end{array}
\right.
\end{equation*}
Clearly $A_{q}(x)=\sum_{n\le x}\chi_{q}(n)$.
Note that $\chi_{q}(n)$ is a completely multiplicative function
in $n$, i.e., $\chi_{q}(nm)=\chi_{q}(n)\chi_{q}(m)$ for all natural
numbers $n$ and $m$. This observation reduces the study of
$\chi_{q}(n)$ to that of $\chi_{q}(r)$ with $r$ a
prime. Using Lemma \ref{Lem a} we infer
the following lemma.

\begin{lem}
We have $\chi_{q}(r)=1$ if and only if $r=p$ or ord$_{r}(q)\not\equiv 2 \; (\mbox{mod } 4)$ in
case $r\ne p$.
\end{lem}

\noindent This result allows one to count the number of primes $r\le x$ such
that $(r,S(q))=1$. Recall that Li$(x)$, the logarithmic integral, is defined
as $\int_2^x{dt/\log t}$.

\begin{lem}
\label{lemma3}
Write $q=p^{t}$. Let $\lambda=\nu_{2}(t)$.
\begin{enumerate}
\item We have
\begin{equation}
\label{noname}
\sum_{r\le x,\;(r,S(q))=1}1=\sum_{r\le x}\chi_{q}(r)=
(1-\delta(q)){\rm Li}(x)+O_{q}\left({x(\log \log x)^4\over \log^3 x}\right).
\end{equation}
\item Assuming GRH the estimate {\rm (\ref{noname})} holds with
error term
$O_{q}(\sqrt{x}\log^{2} x)$, where the index $q$ indicates that
the implied constant depends at most on $q$.
\end{enumerate}
\end{lem}

\pf  1.) The number of primes $r\le x$ such that
ord$_{r}(q)\equiv 2 \; (\mbox{mod } 4)$ is counted in Theorem 2 of \cite{moree}.
On invoking
the Prime Number Theorem in the form $\pi(x)={\rm Li}(x)+O(x\log^{-3}x)$, the
proof of part 1 is then completed.

2.) The proof of this part follows from Theorem 3 of \cite{Moree2}
together with the well-known result (von Koch, 1901) that the Riemann Hypothesis
is equivalent with $\pi(x)={\rm Li}(x)+O(\sqrt{x}\log x)$.
\hfill \Q{\Box}

\vskip .1in

\noindent We are now ready to prove Theorem \ref{main}.

\vskip .1in

\noindent {\bf Proof of Theorem \ref{main}:} 1.) This is a consequence
of part 1 of Lemma \ref{lemma3}, Theorem 4 of \cite{moree} and the fact that
$\chi_q(n)$ is multiplicative in $n$.

2.) By part 2 of Lemma \ref{lemma3} we have
$\sum_{r\le x}\chi_{q}(r)=
(1-\delta(p^{t})){\rm Li}(x)+O_{q}\left({x\log^{-1-v}x}\right)$.
Now invoke Theorem 6 of \cite{MC} with $f(n)=\chi_{q}(n)$.    \hfill \Q{\Box}

\subsection{Counting duadic codes}

Theorem \ref{Thm 6.3.2} allows one to study how many duadic codes of
length $n\le x$ (with $(n,q)=1$) over $\FF_{q}$ exist as $x$ gets large.
We let $D_{q}(x)$ be the associated counting function.
Indeed, we will study the more general function $D_{a}(x)$ which is defined
similarly, but where $a$ is an arbitrary integer. The trivial case arises
when $a$ is a square and thus we assume henceforth that $a$ is not a square.

\indent At first glance it seems that
\begin{equation*}
D_{a}(x)={1\over 2}\sum_{n\le x,\;(n,a)=1}\left(1+({a\over n})\right),
\end{equation*}
with $(a/n)$ the Jacobi symbol. However, it is not true that
$(a/n)=1$ if and only if $a$ is a square modulo $n$, e.g., $(2/15)=(2/3)(2/5)=(-1)(-1)=1$,
but $2$ is not a square modulo $15$. It is possible, however, to develop
a criterium for $a$ to be a square modulo $n$ in terms of Legendre symbols.
To this effect first note that if $a$ is a square modulo $n$, then $a$
must be a square modulo all prime powers in the factorisation of $n$.
This is a consequence of the following lemma.

\begin{lem}
\label{dima}
Let $n$ and $m$ be coprime integers.
Then $a$ is a square modulo $mn$ if and only if it is a square
modulo $m$ and a square modulo $n$.
\end{lem}

\pf  By the Chinese Remainder Theorem
$\mathbb Z/mn\mathbb Z \rightarrow \mathbb Z/m \oplus \mathbb Z/n$ is an
isomorphism of rings and hence $a$ is a square in the
ring on the left if and only if $a$ is is
square in the ring on the right. Now note
that the multiplication in the second ring
is coordinatewise. \hfill \Q{\Box}

\vskip .1in

\noindent It is a well-known result from elementary number theory that
if $p$ is an odd prime and if $x^2\equiv a \; (\mbox{mod } p)$ is solvable, so
is $x^2\equiv a \; (\mbox{mod } p^e)$ for all $e\ge 1$, see e.g.
\cite[Proposition 4.2.3]{IR}. Using this observation together with
Lemma \ref{dima} one arrives at the following criterium for $a$ to be
a square modulo $n$.

\begin{lem}
\label{karakteristiek}
Let $a$ and $n$ be coprime integers.
Put
\begin{equation*}
g_{a}(n)=\prod_{p|n}\left({1+({a\over p})\over 2}\right).
\end{equation*}
Let $e=\nu_{2}(n)$.
Put 
\begin{equation*}
f_{a}(n)=\left\{
\begin{array}{ll} 
0 & \mbox{ if } a\equiv 3 \; (\mbox{mod } 4) \mbox{ and } e\ge 2 ; \\
0 & \mbox{ if } a\equiv 5 \; (\mbox{mod } 8) \mbox{ and } e\ge 3; \\
g_{a}(n) & \mbox{ otherwise}.
\end{array}
\right.
\end{equation*}
Then
\begin{equation*}
f_{a}(n)=\left\{ 
\begin{array}{ll}
1 & \mbox{ if } a \mbox{ is a square modulo } n ; \\
0 & \mbox{ otherwise}.
\end{array}
\right.
\end{equation*}

\end{lem}

\noindent By Lemma \ref{karakteristiek} we have that $D_{a}(x)=\sum_{n\le x,\;(n,a)=1}f_{a}(n)$.
Note that $g_{a}(n)$ is a multiplicative function, but that
$f_{a}(n)$ is a multiplicative function only on the odd integers $n$ (generically).
For this reason let us first consider
\begin{equation*}
G_{a}(x):=\sum_{n\le x,\;(n,a)=1}g_{a}(n).
\end{equation*}
\noindent As a consequence of the law of quadratic reciprocity, the primes $p$
for which $g_{a}(p)=1$ are precisely the primes $p$ in certain arithmetic progressions
with modulus dividing $4q$. On using the prime number theorem for arithmetic progressions
one then infers that for every $v>0$ the following estimate holds true:
\begin{equation}
\sum_{p\le x}g_{a}(p)={1\over 2}{\rm Li}(x)+O_{q}({x\over \log^{v}x}),
\end{equation}
\noindent On using this one sees that the conditions of Theorem 6 of \cite{MC} are
satisfied and this yields the truth of the following assertion.

\begin{lem}
\label{hulp}
Let $\epsilon>0$ and $v\ge 1$ be
arbitrary. Suppose that $a$ is not a square. We have
\begin{equation*}
G_{a}(x)=
\sum_{0\le j <v}{d_{j}x\over \log^{1/2+j}x}+O_{\epsilon,q}\left({x\over
\log^{1/2+v-\epsilon}x}\right),
\end{equation*}
where the implied constant depends at most on $\epsilon$ and $a$, and
$d_{0} \;(>0),\ldots,d_{v}$ are constants that depend at most on $a$.
\end{lem}

\noindent Now it is straightforward to derive an asymptotic for $D_{a}(x)$. Using
Lemma \ref{karakteristiek} one infers that
\begin{equation}
\label{splittje}
 D_{a}(x)=\left\{
\begin{array}{ll}
G_{2a}(x)+G_{2a}({x/2})                  & \mbox{ if } a\equiv 3 \; (\mbox{mod } 4);\\
G_{2a}(x)+G_{2a}({x/2})+G_{2a}({x/4}) & \mbox{ if } a\equiv 5 \; (\mbox{mod } 8);\\
G_{a}(x)                                   & \mbox{ otherwise}.
\end{array}
\right.
\end{equation}

From this and Lemma \ref{hulp} it then follows that we have the following asymptotic for $D_{a}(x)$.
\begin{thm}
\label{ana}
Let $\epsilon>0$ and $v\ge 1$ be
arbitrary. Suppose that $a$ is not a square. We have
\begin{equation*}
D_{a}(x)=
\sum_{0\le j <v}{e_{j}x\over \log^{1/2+j}x}+O_{\epsilon,q}\left({x\over
\log^{1/2+v-\epsilon}x}\right),
\end{equation*}
where the implied constant depends at most on $\epsilon$ and $a$, and
$e_0 \;(>0),\ldots,e_v$ are constants that depend at most on $a$.
\end{thm}

\noindent In particular we have, as $x$ tends to infinity,
\begin{equation*}
D_{a}(x)\sim D_{a}{x\over \sqrt{\log x}} \mbox{ {\rm and} } G_{a}(x)\sim G_{a}{x\over \sqrt{\log x}}\: ,
\end{equation*}
where $D_{a}$ and $G_{a}$ are positive constants. We now consider
the explicit evaluation of these
constants. Note that by (\ref{splittje}) it suffices
to find an explicit formula for the constant $G_{a}$.

In case $a=D$ is a negative discriminant of a binary quadratic form this constant
can be easily computed using results from the analytic theory of binary quadratic
forms. We say an integer $D$ is a discriminant if it arises as the discriminant
of a binary quadratic form. This implies that either $4|D$ or $D\equiv 1 \; (\mbox{mod } 4)$.
On the other hand, it can be shown that any number $D$ satisfying
$4|D$ or $D\equiv 1 \; (\mbox{mod } 4)$ arises as the discriminant of a binary quadratic
form. Now let $D$ be a discriminant and $\xi_{D}$ be the multiplicative function defined as
follows:
\begin{equation*}
\xi_{D}(p^e)=\left\{ 
\begin{array}{ll}
1 & \mbox{ if } ({D\over p})=1 ;\\
1 & \mbox{ if } ({D\over p})=-1 \mbox{ and } 2|e;\\
0 & \mbox{ otherwise}.
\end{array}
\right.
\end{equation*}
\noindent Let $n$ be any integer coprime to $D$. Then $\xi_{D}(n)=1$ if and only if $n$ is represented
by some primitive positive integral binary quadratic form of discriminant $D$. Let
$B_{D}(x)$ denote the number of positive integers $n\le x$ which are coprime to $D$
and which are represented by some primitive integral form of discriminant $D\le -3$.
Note that $B_{D}(x)=\sum_{n\le x}\xi_{D}(n)$. It was proved by James \cite{J} that
\begin{equation*}
B_{D}(x)=J(D){x\over \sqrt{\log x}}+O\left({x\over \log x}\right),
\end{equation*}
where $J(D)$ is the positive constant given by
\begin{equation}
\label{jd}
\pi J(D)^2={\varphi(|D|)\over |D|}L(1,\chi_{D})\prod_{({D\over p})=-1}
{1\over 1-{1\over p^2}},
\end{equation}
and $p$ runs over all primes such that $(D/p)=-1$.
(Recall that the Dirichlet L-series $L(s,\chi_{D})$ is
defined by $L(s,\chi_{D})=\sum_{n=1}^{\infty}\chi_{D}(n)n^{-s}$.)
Since the behaviour of $\xi_{D}$
is so similar to that of $f_{D}$, James' result can in fact be used to determine
the asymptotic behaviour of $G_{D}(x)$ for negative
discriminants $D$ and, in particular, to determine
$G_{D}$. Using a classical result of Wirsing, see
e.g. Theorem 3 of \cite{MC}, one infers that
\begin{equation*}
{G_{D}(x)\over B_{D}(x)}\sim \prod_{p\le x\atop ({D\over p})=-1}\left(1-{1\over p^2}\right).
\end{equation*}
{}From this and the identity (\ref{jd}) it follows that $G_{D}$ is the positive solution of
\begin{equation}
\label{deetje}
\pi G_{D}^2={\varphi(|D|)\over |D|}L(1,\chi_{D})\prod_{({D\over p})=-1}
\left(1-{1\over p^2}\right).
\end{equation}

\noindent For more details on $B_{D}(x)$ and
related counting functions the reader is referred to a
paper (in preparation) by Moree and Osburn \cite{MO}. In \cite{MO} it is also
pointed out that $B_{D}(x)$ in fact satisfies an asymptotic result similar to
the one given for
$D_{a}(x)$ in Theorem \ref{ana}.

The fact that the characteristic functions $\xi_{d}$ and $f_{D}$ are so closely
connected, can be exploited to give a criterium for the existence of duadic codes
in terms of representability by quadratic forms.

\begin{lem}
Let $q$ be an odd prime power, say $q=p_{1}^{e}$ with $p_{1} \equiv 3 \; (\mbox{mod } 4)$. Let $n$ be an 
odd squarefree integer satisfying $(n,q)=1$ and suppose, moreover, that
$n$ can be written as a sum of two integer squares. A duadic code of length $n$ over
$\FF_{q}$ exists if and only if $n$ can be represented by
some primitive positive integral binary quadratic form of discriminant $-p_{1}$
\end{lem}

\pf By assumption $-p_{1} \equiv 1 \; (\mbox{mod } 4)$ and hence is
a discriminant. The assumption that $n$ is odd and squarefree ensures that
$\xi_{-p_{1}}(n)=f_{-p_{1}}(n)=f_{-p_{1}^{e}}(n)$. The assumption that $n$ can be
represented as a sum of two squares, together with the assumption that $n$ is squarefree
ensures that $n$ is a product of primes $p$ satisfying $p\equiv 1 \; (\mbox{mod } 4)$.
For every prime $p\equiv 1 \; (\mbox{mod } 4)$ we have $(-p_{1}^{e}/p)=(p_{1}^{e}/p)$. It
thus follows that $\xi_{-p_{1}}(n)=f_{-p_{1}^{e}}(n)=f_{p_{1}^{e}}(n)$. The result then
follows on invoking Theorem \ref{Thm 6.3.2}, Lemma \ref{karakteristiek} and
the fact that, for $(n,D)=1$, $\xi_{D}(n)=1$ if and only if $n$ is represented by some
primitive positive integral binary quadratic form of discriminant $D$. \hfill \Q{\Box}\\

\noindent It remains,  however, to determine $G_{a}$ for a general number $a$. It
is well-known from Tauberian theory that one has
\begin{equation*}
G_{a}={1\over \Gamma(1/2)}\lim_{s\downarrow 1}\sqrt{s-1}F(s),
\end{equation*}
where $F(s)=\sum_{n=1}^{\infty}g_{a}(n)n^{-s}$. An easy computation shows that
\begin{equation*}
(s-1)F(s)^2=(s-1)\zeta(s)
{\varphi(|a|)\over |a|}
L(s,\chi_{a})\prod_{({a\over p})=-1}
\left(1-{1\over p^2}\right).
\end{equation*}
On using that the Riemann zeta-function $\zeta(s)$ has a simple pole at $s=1$ of
residue $1$, one obtains that
\begin{equation*}
\pi G_{a}^2={\varphi(|a|)\over |a|}L(1,\chi_{a})\prod_{({a\over p})=-1}
\left(1-{1\over p^2}\right).
\end{equation*}
Notice that equation (\ref{deetje}) is a special case of this.

\vskip .5in
\noindent {\bf Acknowledgements}. The first author gratefully acknowledges financial support from the 
University of the Philippines and from the Philippine Council for Advanced Science and Technology Research 
and Development through the Department of Science and Technology.

\indent The second author likes to thank I. Shparlinski for suggesting Lemma A.5 and D. Gurevich for some 
helpful remarks.

\indent The three authors would like to thank the anonymous referees for their helpful suggestions that 
greatly improved the presentation of the material.

\end{document}